\documentclass[a4paper]{amsart}

\usepackage[all]{xy}
\usepackage{amsmath}
\usepackage{amsfonts,amssymb}
\numberwithin{equation}{section}
\newtheorem{satz}{Proposition}
\newtheorem{Satz}[satz]{Theorem}
\newtheorem{definition}{Definition}
\newtheorem{lemma}[satz]{Lemma}

\newtheorem{kor}[satz]{Corollary}

\theoremstyle{definition}
\newtheorem*{remark}{Remark}
\newtheorem*{notation}{Notation}

\newcommand{\tensor}{\otimes}

\newcommand{\un}[1]{\ensuremath{\protect\underline{#1}}}

\def\1halb{\frac{1}{2}}

\def\tto{\twoheadrightarrow}

\DeclareMathOperator{\Spec}{Spec}

\DeclareMathOperator{\Aut}{Aut}
\DeclareMathOperator{\End}{End}

\DeclareMathOperator{\Char}{char}

\DeclareMathOperator{\Bun}{Bun}

\def\sxymat{\xymatrix@C=1.5ex@R=0.8ex}
\def\grp{$\xymatrix{ R\times_{X}R  \ar[r]^-{\mu} & R \ar@<1ex>[r]^-{s}\ar@<-1ex>[r]_-{t} & X}$}
\def\dar{\ar@<-0.5ex>[r]\ar@<0.5ex>[r]}
\def\tar{\ar[r]\ar@<1ex>[r]\ar@<-1ex>[r]}
\newcommand{\dmap}[2]{\ar@<-0.5ex>[r]_-{#2}\ar@<0.5ex>[r]^-{#1}}
\newcommand{\dotarrow}[2]{\xymatrix{{#1}\ar@{..>}[r]&{#2}}} 
\def\cart{\ar@{}[dr]|{\square}}

\def\cG{\mathcal{G}}
\def\cH{\mathcal{H}}

\def\cL{\ensuremath{\mathcal{L}}}

\def\cO{\mathcal{O}}
\def\cP{\mathcal{P}}
\def\cQ{\mathcal{Q}}

\def\cg{\mathfrak{g}}
\def\cm{\mathfrak{m}}

\def\cp{\mathfrak{p}}
\def\cq{\mathfrak{q}}

\def\cu{\mathfrak{u}}

\def\bA{{\mathbb A}}

\def\bG{{\mathbb G}}

\def\bZ{{\mathbb Z}}

\def\bR{{\mathbf R}}

\begin{document}
\title[Langton's theorem]{Semistable reduction for $G-$bundles on curves}
\author{Jochen Heinloth}
\address{Fachbereich Mathematik der Universit\"at Duisburg-Essen, D-45117 Essen, Germany}
\email{jochen.heinloth@uni-due.de}
\begin{abstract}
We prove a  semistable reduction theorem for principal bundles on curves in almost arbitrary characteristics. For exceptional groups we need some small explicit restrictions on the characteristic.
\end{abstract}
\maketitle

\section{Introduction}

We want to prove a semistable reduction theorem for principal $G$-bundles over a smooth projective curve.

\begin{Satz}\label{semistab}
Let $C_{/k}$ be a smooth projective curve over an algebraically closed field $k$ and let $G$ be a semisimple group over $k$. 

Assume either that Behrend's conjecture holds for $G$ or assume that $\Char(k)\neq 2$ if $G$ contains factors of type $B,C,D$, $\Char(k)>7$ if $G$ contains a factor of type $G_2$, $\Char(k)>19$ if $G$ contains a factor of type $F_4$ or $E_6$, $\Char(k)>31$ if $G$ contains a factor of type $E_7$ and $\Char(k)>53$ if $G$ contains a factor of type $E_8$.

Let $R$ be a complete discrete valuation ring over $k$ with fraction field $K$.  Then for any semistable principal $G$-bundle $\cG_K$ on $C\times_k \Spec K$ there exists a finite extension $R^\prime$ of $R$ such that the pullback of $\cG_K$ extends to a semistable $G$-bundle on $C\times_k \Spec R^\prime$.
\end{Satz}
Recall that Behrend's conjecture (\cite{BehrendInstability}, 7.6) says that the reduction of a principal bundle to the maximal instability parabolic has no infinitesimal deformations.
This is only known in characteristic $0$ and for characteristics bigger than the height $h(G)$, by Biswas and Holla \cite{BiswasHolla}. For exceptional groups we need to apply this result in the last step of our proof, which causes the bounds on the characteristics. In case $G$ is of type $A$, our proof also works in characteristic $2$, as in Langton's theorem \cite{Langton}.

If the characteristic of $k$ is $0$, then the above theorem has been shown by Faltings \cite{FaltingsBundles} and Balaji and Seshadri \cite{BalajiSeshadri}. This was generalized to large characteristics by Balaji and Parameswaran \cite{Balaji}.
Our motivation to look for another approach to this theorem came from the recent work of Gomez, Langer, Schmitt and Sols \cite{GLSS}, who managed to construct quasi--projective coarse moduli spaces for semistable principal bundles in arbitrary characteristic. The above theorem shows that in the case of smooth projective curves their moduli spaces are projective in almost any characteristic.

Our strategy is similar to Langton's proof \cite{Langton} of the corresponding theorem for vector bundles:
\begin{enumerate}
\item Take any extension $\cG_R$ of $\cG_K$ to $C\times_k \Spec R$. This is possible -- at least after a finite extension of $R$ -- because the affine Gra\ss mannian is (ind-)projective.
\item If the special fiber $\cG_k$ of the principal bundle is not semistable, then we can apply Behrend's theorem to find a canonical reduction of $\cG_k$ to a parabolic subgroup $P\subset G$. Since $\cG_K$ is semistable there is a maximal $n$ such that this canonical reduction to $P$ can be lifted to $R/\cm^n$, where $\cm$ is the maximal ideal of $R$.

Choose a cocycle for $\cG_R$, such that it reduces to a cocycle defining the canonical reduction modulo $\cm^n$. (We even have to be a little more careful in choosing the cocycle, see Proposition \ref{cocycle}.) Conjugating the cocycle for $\cG_R$ with a $K$-valued element of the center of the Levi subgroup $L$ of $P$ we find a modification $\cG^\prime_R$ of $\cG_R$ such that the special fiber of $\cG^\prime_R$ has a reduction to the opposite parabolic $P^-$,  the associated $L$-bundle is isomorphic to the $L$-bundle defined by the canonical reduction of the special fiber of $\cG$ and finally this $P^-$-bundle has no reduction to $L$.

\item Deduce that the special fiber of $\cG^\prime_R$ is less unstable than the special fiber of $\cG_R$. Here ``less unstable'' means that the bundle lies in a component of the Harder-Narasimhan-stratification of the moduli stack of $G$-Bundles, which contains $\cG_R$ in its closure. This holds, since by our construction $\cG^\prime_k$ has a reduction to a parabolic $P^-$, such that the associated $L$-bundle is the same as the one obtained form the canonical reduction of $\cG_k$. We show that the deformation which deforms $\cG^\prime_k$ into the bundle induced from the associated $L$-bundle, is transversal to the Harder-Narasimhan-stratum of $\cG_k$. Only here we need our extra condition on $G$, because we have to avoid higher order deformations of the canonical parabolic reduction.
\end{enumerate}

Let us compare this to Langton's original argument. He performs a modification of the vector bundle, such that the new bundle has a quotient which is isomorphic to the maximal destabilizing subsheaf of the old bundle. This corresponds to the reduction to $P^-$ in our situation. For vector bundles this new bundle is less unstable, unless the above quotient is a direct summand, which can happen only if the destabilizing subsheaf of the original bundle can be lifted to $C\times \Spec(R/\cm^2)$. In our construction we therefore use a modification of higher order in order to obtain a bundle for which the quotient defines a nontrivial extension (condition (4) in Proposition \ref{cocycle}). This property allows us to prove that the new bundle is always less unstable.

A second difference to Langton's algorithm is that we use the affine Gra\ss mannian and cocycles instead of the Bruhat-Tits building. This has the advantage that we can stay within the category of $G$-bundles and do not have to consider more complicated families of group schemes over $C$.

\noindent {\bf Acknowledgments.} I would like to thank Alexander Schmitt for pointing out the problem of semistable reduction for principal bundles to me and for many discussions and comments. I would like to thank G. Harder, U. Stuhler and V. Balaji for their comments. I would like to thank A. Langer and Ngo Tuan Doc for pointing out a mistake in a previous version of this article.

\begin{notation} Throughout the article we will fix $C_{/k}$ a smooth projective curve over an algebraically closed field $k$, $G$ a semisimple group scheme over $k$. $R$ denotes a complete discrete valuation ring over $k$ with residue field $k$ and quotient field $K$. Finally, let $\pi$ be a local parameter for the maximal ideal $\cm\subset R$. We will often need to replace $R$ by a finite extension, which we will then denote by $R$ again.  

For a $k-$algebra $A$ we will denote the base-change of any object over $k$ to $A$ by a lower index $A$, e.g. $C_A:=C\times_k\Spec A$.

Given an affine algebraic group $H$ which acts on a scheme $S$ by $\rho:H\to \Aut(S)$ and an $H$-bundle $\cH$ on $C$, we will denote by $\cH\times^{H,\rho} S:=\cH \times S/H$ the associated bundle with fiber $S$.

For an algebraic groups $G,P,Q,U$ we will denote the corresponding Lie algebras by $\cg,\cp,\cq,\cu$.

The group of $\ell$-th roots of unity is denoted by $\mu_\ell$.
\end{notation}

\begin{remark}
Since the algorithm will terminate after a finite number of steps, the assumption that the residue field of $R$ is $k$ and therefore algebraically closed is not essential,  we could as well start with an arbitrary complete discrete valuation ring, and take the finite extension needed in each step to find a canonical parabolic subgroup in the special fiber.
\end{remark}

\section{Stability of $G$-bundles}
Let us recall the concept of stability for $G$-bundles.
\begin{definition}(Ramanathan) A principal $G$-bundle $\cG$ on a curve $C$ is called {\em semi\-stable} if for all parabolic subgroups $P\subset G$, all reductions of $\cG$ to $\cP$ and all dominant characters $\alpha: P\to \bG_m$ we have:
$$ \deg(\cP\times^{P,\alpha} \bA^1) \leq 0.$$
\end{definition}
Here $\cP\times^{P,\alpha} \bA^1$ is the line bundle obtained from the action of $P$ on $\bA^1$ defined by $\alpha$.

Furthermore Behrend shows \cite{BehrendInstability} (Proposition 8.2/Theorem 7.3) that every $G$-bundle $\cG$ over $C$ has a canonical reduction $\cP$ to a parabolic subgroup $P\subset G$, which is characterized by:
\begin{enumerate}
\item Let $R_u(P)$ be the unipotent radical of $P$. Then $\cP/R_u(P)$ is a semi--stable $P/R_u(P)$ bundle.
\item  For all dominant characters $\alpha:P\to\bG_m$ we have $\deg(\cP\times^{P,\alpha} \bA^1)>0$.
\end{enumerate}

The {\em type $t(\cP)$ of the canonical reduction} of $\cG$ is given by the type of the parabolic subgroup $P$ together with the degree of $\cP$, which is defined as 
$$ \un{\deg}: X^*(P) \to \bZ$$
$$  \alpha \mapsto \deg(\cP\times^{P,\alpha} \bA^1).$$

The type of the canonical reduction induces a stratification of the moduli stack of principal $G$-bundles on $C$, the Harder-Narasimhan-stratification. It has become tradition to call this decomposition into constructible subsets stratification even if it is not obvious that the closure of a stratum is a union of strata. However, the dimension of a stratum depends only on the degree of instability, which is upper continuous. Further, the uniqueness of the canonical reduction implies that a specialization of a bundle of given type of instability cannot intersect a different stratum of the same degree of instability. Therefore the closure of a stratum can only intersect strata of lower dimension.
We will denote by $\Bun_G$ the moduli stack of $G$-bundles, similarly $\Bun_P^{\un{d}}$ is the moduli stack of $P$-bundles of degree $\un{d}$. The substack of semistable bundles will be denoted by $\Bun_P^{\un{d},ss}$. 

\begin{lemma}\label{adjoint}
It is sufficient to prove our theorem for adjoint groups.
\end{lemma}
\begin{proof}
From the definition of semistability we see that a $G$-bundle is semistable if and only if the induced $G/Z(G)$-bundle (where $Z(G)$ is the center of $G$) is semistable.
Moreover the map $\Bun_G\to \Bun_{G/Z(G)}$ is finite, because the obstruction to lift a $G/Z(G)$-bundle to a $G$-bundle is given by a class in $H^2(C,Z(L))$, so the image of $\Bun_G$ in $\Bun_{G/Z(G)}$ is closed and the fibres are a torsor for $\Bun_{Z(G)}$, which is finite.
\end{proof}

\begin{remark} Since a semisimple adjoint group $G_{/k}$ is a product of simple groups (\cite{SGA3} Expos\'e XXIV, Prop. 5.9), it would be sufficient to prove our theorem for simple adjoint groups. We will use this only to simplify the last step of our proof for groups containing an exceptional factor.
\end{remark}
\section{First step: Find an arbitrary extension of $\cG_K$} 
Recall the following theorems: Let $S$ be any scheme over $k$.

\begin{Satz}(Drinfeld--Simpson, Harder)
Let $G/k$ be semisimple and $U\subset C$ an affine open subset. Then the restriction of any $G$-bundle on $C\times S$ to $U\times S$ becomes trivial after a suitable faithfully flat base change $S^\prime \to S$, which is locally of finite presentation. If $G$ is simply connected, then $S^\prime \to S$ can be chosen to be \'etale.
\end{Satz}
See \cite{DrinfeldSimpson}, or \cite{HarderDedekind} for $S=\Spec k$.

\begin{Satz}(Affine Gra\ss mannian)
There is an ind-projective scheme $G((t))/G[[t]]$, which parameterizes $G-$bundles on $C$ together with a trivialization on $C-p$, where $p$ is a closed point on $C$.
\end{Satz}
I am not sure who was the first to find this result. A proof general enough to include the above formulation can be found in \cite{FaltingsLoopGroups} (Corollary 3 and the remark following the proof of this corollary), see also \cite{BD} 2.3.4 and 4.5.1.

These theorems show that we can extend $G-$bundles: Let $\cG_K$ be  a $G-$bundle on $C_K$. Choose a point $p\in C(k)$, then by the first theorem we can replace $K$ by a finite extension (and $R$ by its integral closure in this extension), such that $\cG_K|_{(C-p)\times \Spec(K)}$ becomes trivial. Any trivialization defines a map $\Spec(K)\to G((t))/G[[t]]$. Since the affine Gra\ss mannian is ind-projective this defines a map $\Spec(R)\to G((t))/G[[t]]$ and thus a $G$-bundle $\cG_{R}$ on $C_R$, which extends $\cG_K$.

\begin{remark} For any linear algebraic group $H$ over $k$ there is an ind-scheme $H((t))$, called loop group, which is defined by $H((t))(R)=H(R((t)))$ for all $k$-algebras $R$ (e.g. \cite{FaltingsLoopGroups}). In the following we will often want to lift points $G((t))/G[[t]](R)$ to $G((t))(R)$. This is possible -- after possibly extending the coefficients -- by choosing a trivialization of the $G$-bundle on the formal completion at $p$, the glueing datum for the $G$-bundle is then a point of $G((t))$.

Recall that conversely, any element of $g\in G((t))(R)$ defines a glueing datum for a $G$-bundle on $C$. For noetherian $R$ this comes from flat descent, in general this was shown by Beauville-Laszlo \cite{BeauvilleLaszlo}. From this description one also sees that if we take $g_{C-p} \in G(R[C-p])$ and $g_p\in G(R[[t]])$ then the bundles given by $g$ and $g_{C-p}gg_p$ are isomorphic. 
\end{remark}

\section{Second step: Find a modification of $\cG_R$}

Assume that the special fiber $\cG_k$ of $\cG_R$ is not semistable. Since $k$ is algebraically closed, the bundle $\cG_k$ has a canonical reduction $\cP_{k}$ to a maximal instable parabolic subgroup $P\subset G$ (\cite{BehrendInstability} Thm 7.3, Prop 8.2).

We want to find a cocycle for $\cG_R$, such that its reduction modulo $\pi$ is a cocycle with values in $P\subset G$ which defines $\cP_k$. To do so, first we have to note:
\begin{lemma}
The bundle $\cP_{k}$ is Zariski locally trivial on $C_k$.
\end{lemma}
\begin{proof} 
This follows from a general theorem in SGA 3:
We know that the bundle $\cG_k$ is Zariski locally trivial on $C_k$, because $G$ is semisimple. Therefore we may assume that $\cP_k$ is a reduction to $P$ of the trivial bundle over an open subset $U\subset C$, i.e. a section of $G/P\times U$. 
Now by \cite{SGA3} (Expos\'e XXVI Cor. 5.2) we know that over any (semi-)local base $S$ we have $G/P(S)=G(S)/P(S)$. Therefore there is an open covering  $\cup_{i=1}^n V_i= U$ and $g_i\in G(V_i)$ such that $(\cP_k)|_{V_i} =  g_i.(P \times V_i) \cong P\times V_i$.
\end{proof}

Since $\cP_k$  is Zariski locally trivial, there is a finite set of points $S\subset C(k)$ such that $\cP|_{C-S}$ is trivial and the restrictions of $\cP_k$ to the complete local rings at points in $S$ are trivial as well. To simplify notations we choose local parameters at all points in $S$, i.e. isomorphisms $\cO_{C,p}^{\widehat{\quad}}\cong k[[t]]$. Then every element $\prod_{p\in S} G((t))(R)$  defines a $G$-bundle on $C_R$. With this notation we can conclude:
\begin{kor}\label{cocycle1}
There exists a finite set $S\subset C(k)$ and an element $\overline{p}\in \prod_{p\in S} P((t))(k)$ such that $\cP_k$ is isomorphic to the bundle defined by $\overline{p}$.

Furthermore there exists $g\in \prod_{p\in S} G((t))(R)$ such that $g$ defines the bundle $\cG_R$ and $\overline{p} \equiv g\mod (\pi).$
\end{kor}
\begin{proof}
The first part has been shown above. For the second part we only have to note that since $G$ is smooth, the lifting criterion for smoothness tells us that we can lift any local trivialization of $\cG_k$ to $\cG_R$ (using our simplifying assumption that $R$ is complete). If we start with a local trivialization of the canonical reduction $\cP_k$ and lift this to a trivialization of $\cG_R$ then we get a cocycle $g$ as claimed.
\end{proof}

To find a modification of $\cG_R$ with less unstable special fiber we choose a maximal parabolic $Q\supset P$. The canonical reduction $\cP_k$ defines a reduction $\cQ_k$ of $\cG_k$ to $Q$. As indicated in the introduction we want to choose our cocycle such that it defines a 
lifting of this reduction to $\cG|_{C\times \Spec(R/(\pi^n))}$ for $n$ maximal. Such an $n$ exists, because the generic fiber of $\cG_R$ is semistable, therefore we cannot obtain a lift to $C_R$.
We even have to be a bit more careful and introduce some standard notation for algebraic groups (see \cite{SGA3} Exp. XXII, 1.1 and Exp. XXVI):

We choose an \'epinglage of $G$ adapted to $P$, i.e.: We choose $B\subset P\subset G$, a Borel subgroup of $G$ contained in $P$ and $T\subset B$ a maximal torus, which is split because $k$ was assumed to be algebraically closed. For all roots $\alpha$, we denote by $U_\alpha\subset G$ 
the root-subgroup of $G$ (it is isomorphic to $\bG_a$ and conjugation by $T$ acts on $U_\alpha$ via the character $\alpha$).

Let $\Delta\subset X^*(T)$ be the set of roots of $T$, $\Delta^+$ the set of positive roots, $I\subset \Delta^+$ the set of positive, simple roots determined by $B$. We denote by $\Delta_{Q}\subset \Delta$ the roots of $Q$ and $I_Q:= I\cap \Delta_{Q}$ the  simple roots contained in $Q$ and let $\beta\in I-I_Q$ be the unique positive simple root, not contained in $I_Q$.

Let $L\subset Q$ be a Levi subgroup. Then $Z(L)^\circ:=(\cap_{\alpha\in I_Q} \ker(\alpha))^\circ$ is the connected component of the center $Z(L)$ of $L$. Since the positive simple roots form a basis of $X^*(T)\tensor_\bZ \bR$ we have surjections with finite kernels:
$$\xymatrix{
T \ar@{->>}[d]^-{\prod \alpha}& Z(L)^\circ\cong \bG_m \ar@{->>}[d]^{\lambda:=\beta|_{Z(L)^\circ}}\ar@{_(->}[l]\\
 \prod_{\alpha\in I} \bG_m & \bG_m\ar@{_(->}[l]_-{id_\beta}\\
}
$$
Recall that the unipotent radical $U$ of $Q$ has a filtration $U=U_1 \supset U_2 \supset \dots \supset U_h$ by normal subgroups, where $U_i:= \prod_{\alpha\in \Delta_{\beta,i}} U_\alpha$, where $\Delta_{\beta,i} = \{ \alpha=\sum_{k\in I_Q} n_k \alpha_i + m \beta \in \Delta^+|  m \geq i\}$. Thus this filtration can also be obtained by decomposing $U$ into eigenspaces for the action of $Z(L)^\circ$ which acts by conjugation on $U$. 

Similarly for the opposite parabolic $Q^-$ of $Q$ (defined by $-\Delta_P$), we have a filtration $U^-=U^-_1\supset U_2^- \supset \dots\supset U_h^- $. Finally denote by $LU_i^-:=L\ltimes U_i^- \subset Q^-$ the subgroup generated by $L$ and $U_i^-$. 

\begin{remark} Note that the length of the above filtration is bounded by the largest coefficient occurring when writing the roots as linear combination of simple roots. In particular for all groups of type $B,C,D$ only coefficients $\leq 2$ occur (\cite{BourbakiLie}), so $U_3=(1)$ in these cases.
For groups of type $A$ we have $U_2=(1)$.
\end{remark}

To define our modification of $\cG_R$ we will have to extend $R$ to $R[\pi^{1/N}]$, where $N$ is chosen such that for all $i\in\{1,2,\dots,h\}$ there exists $k_i$ such that $(\lambda(\pi^{k_i/N})))^i=\pi$.

Let $\cL_k:=\cQ_k\times^Q L=\cQ_k/U$ be the $L$-bundle obtained from $\cQ_k$. Then any $LU_i^-$-bundle $\cQ^-_{i}$ on $C$  with $\cQ^-_i/U^-_i \cong \cL$ defines an element in $H^1(C, ^\cL U_i^-)$, where $^\cL U_i^-:=\cL \times^L U_i^- $ is the group scheme over $C$ defined by the conjugation action of $L$ on $U_i$.

\begin{satz}\label{cocycle}
There is a cocycle $g\in \prod G((t))(R)$ for $\cG_R$ and $z=\pi^{\ell/N}\in \bG_m(K)=Z(L)^\circ(K)$ such that the following holds:
\begin{enumerate}
\item $g \mod \pi \in \prod_{p\in S} P((t))(k)$ and this cocycle defines the canonical reduction of $\cG_k$ to $P$.
\item $g^\prime:= zgz^{-1} \in \prod_{p\in S} G((t))(R[\pi^{1/N}])$.
\item  Denote $\overline{g}^\prime:= g^\prime \mod \pi^{1/N}$. Then $\overline{g} ^\prime\in \prod_{p\in S} Q^-((t))(k)$.
\item Let $i$ be maximal, such that $\overline{g}^\prime \in \prod_{p\in S} (LU_i^-((t))(k)$. Then the class $[\overline{g}^\prime]\in H^1(C,^{\cL}U^-_i/^\cL U^-_{i+1})$ is non-zero.
\end{enumerate}
\end{satz}

\begin{proof}
We will inductively modify the cocycle $g$ we found in Corollary \ref{cocycle1}.
First note that each component of $g$ is an $R((t))$-valued point of $G$ contained in the open subset $U^- \times L \times U \subset G$, since $g\mod \pi \in \prod_{p\in S} P((t))(k)$.   Therefore we can decompose $g = v\cdot  l \cdot u$ with $v\in \prod_{p\in S} U^-(R((t)))$, $l\in\prod_{p\in S} L(R((t)))$, $u\in \prod_{p\in S} U(R((t)))$.

Let $\ell$ be  maximal such that for $z:=\pi^{\ell/N}$ we have that  $g^\prime := z g z^{-1}$ is an $R[\pi^{1/N}]((t))$-valued point of $\prod_{p\in S} G$. Then $\ell>0$ since $v\mod \pi = 1$, by property (1) of $g$. Furthermore, if $n$ is the maximal integer such that we can lift of the canonical reduction of $\cG_k$ to $Q$ to $R/\pi^n$, then $\ell\leq nN$, because otherwise (2) implies that $g\mod \pi^{n+1} \in \prod_{p\in S} Q((t))(R/(\pi^{n+1}))$, contradicting the maximality of $n$. Therefore we will show, that if $g$ does not satisfy the conditions of the lemma, then we can find another cocycle for which we can increase $\ell$.

Write $\overline{v}^\prime:= z v z^{-1} \mod \pi^{1/N}$ and $\overline{l}:= l \mod \pi^{1/N}$. Then $\overline{g}^\prime = \overline{v}^\prime \overline{l}$.

Let $i$ be maximal such that all components of $\overline{v}^\prime$ lie in $U_i^-$. In particular $U_i^-\neq\{ 1\}$, otherwise we could increase $\ell$.
If the class $[\overline{v}^\prime \overline{l}]\in H^1(C,^{\cL}U_i^-/^\cL U_{i+1}^-)$ is zero, then the cocycle $\overline{v}^\prime\overline{l}$ is a boundary, i.e. there exist $\overline{v}_{C-S}\in \prod_{p\in S} U_i^-(k[C-S])$, $\overline{v}_S\in \prod_{p\in S} U_i^-(k[[t]])$ such that $ \overline{v}_{C-S} \overline{v}^{\prime } \overline{l} \overline{v}_S = \overline{v}_{i+1} \overline{l}$ with $\overline{v}_{i+1}\in \prod_{p\in S} U_{i+1}^-(k((t)))$.

Take any lift of $v_{C-S}\in U_i^-(R[C-S])$ and $v_S\in \prod_{p\in S} U_i^-(R[[t]])$. Then
$ v_{C-S} g^\prime v_S \mod \pi^{1/N} \in \prod LU_{i+1}^-((t))(k)$ and
$$ v_{C-S} g^\prime v_S =z \bigg( \big( z^{-1}v_{C-S} z\big) g \big( z^{-1} v_S z\big)\bigg) z^{-1}.$$

The inner bracket is an $R[\pi^{1/N}]$-valued cocycle for $\cG_{R[\pi^{1/N}]}$. We can replace this by an $R$-valued cocycle as follows: Let $\mu_m:=\ker(\lambda)$. Then conjugation by $z$ acts as multiplication by $z^{-mi}$ on $U_i^-/U_{i+1}^-$. Since $\ell$ and $i$ were chosen maximal $mi\ell/N$ is an integer. Therefore the element $z^{-1}v_{C-S}z$ defines an $R[C-S]$-valued point of $U_{i}^-/U_{i+1}^-$, congruent $1$ mod $\pi^{mi\ell/N}$. The same holds for $z^{-1}v_{S}z$.

Choose $v_{C-S}^\prime \in U_i^-(R[C-S])$ such that $v_{C-S}^\prime  \equiv z^{-1} v_{C-p} z \mod U_{i+1}^-$ and $zv_{C-S}^\prime z^{-1}$ is integral over $R[\pi^{1/N}]$. Similarly choose $v_S^\prime\in \prod_{p\in S} U_i^-(R((t)))$ with the analogous properties.

Define a new cocycle $g_1 := v_{C-S}^\prime g v_S^\prime$. By construction $g_1$ satisfies (1)--(3), but now $z g_1z^{-1} \mod \pi^{1/N}\in \prod_{p\in S} LU_{i+1}^-((t))(k)$, because this holds for $v_{C-S}g^\prime v_S$. If $U_{i+1}^-=\{1\}$, this means that we can find a larger $\ell$ such that 2 and 3 hold. Continuing this process inductively we can find a cocycle $g$ as claimed, since $\ell$ is bounded by $nN$.
\end{proof}

Now we replace $R$ by $R[\pi^{1/N}]$ and let $\cG^\prime_R$  be the $G$-bundle defined by the cocycle $g^\prime$ constructed in the above lemma.

\section{Third Step: Show that the modified bundle $\cG_R^\prime$ is less unstable}

To show that the modified bundle is less unstable we will need to compare the $G$-bundle induced  from a $Q$ (or $Q^-$) bundle to the one induced from the Levi-quotient. Since $Z(L)^\circ\cong \bG_m$ and $\beta\in \Delta_Q^+$, the conjugation action $Z(L)^\circ \times Q \to Q$ given by $(z,q) \mapsto z q z^{-1}$ extends to a map $\widetilde{\lambda}:\bA^1 \times Q \to Q$.

For any $Q$-bundle $\cQ_k$ this induces an action of $Q$ on $\cQ_k \times \bA^1\times Q$, defined as $(x,t,q).q^\prime:= (x.q^\prime,t,\widetilde{\lambda}(t,q^\prime).q)$. The quotient by this action is a $Q$-bundle on $C\times \bA^1$:
$$ \cQ_{\lambda} := \cQ_k \times^{Q}(\bA^1 \times Q).$$
Similarly, for a $Q^-$-bundle $\cQ^-_k$ we use the action of $\bG_m$ on $Q^-$ given by $\lambda^{-1}$ to define:
$$\cQ^-_{\lambda}:= \cQ^-_k \times^{Q^-} (\bA^1 \times Q^-).$$

\begin{lemma}\label{deformation}
There are canonical isomorphisms of $Q$-bundles:
\begin{enumerate}
\item $\cQ_{\lambda}|_{C\times \bG_m} \cong \cQ_k \times \bG_m.$
\item Let $t$ be a coordinate of $\bA^1$ and $\ker(\lambda)\cong \mu_m$, then $\cQ_{\lambda}|_{C\times \Spec k[t]/t^m}\cong \cQ_k/U\times^L Q \times \Spec(k[t]/t^m)$.
\end{enumerate}
\end{lemma}
\begin{proof}
For the first part an isomorphism $\cQ_k \times \bG_m \to \cQ_k \times^Q ( \bG_m \times Q)$
is given by $(x,t) \mapsto (x.\lambda(t),t,1)$.
For the second part we only need to notice that, if $\lambda$ factors through the map $t\mapsto t^l$, then $\widetilde{\lambda}|_{\Spec k[t]/t^l} = i \circ p$, where $p: Q \to Q/U=L$ is the projection and $i:L\to Q$ the inclusion.
\end{proof}
\begin{remark}
If $\cP_k$ is given as the canonical reduction of a $G$-bundle $\cG_k$ to the instability parabolic $P$, and $\cQ_k:=\cP\times^P Q$, then the canonical reduction of $\cG_L:=\cQ_k/U \times^L G$ is given by $\cP_k/(U\cap P)\times^{L\cap P} P$. In particular, $\cG_L$  and $\cG_k$ lie in the same HN-stratum of $\Bun_G$.
\end{remark}
This holds because, by definition of the canonical reduction $\cP_k/R_u(P)$ is a semistable bundle, and the degree of the line bundles obtained from characters of $P$ are constant in the family defined above.
\begin{kor}
The closure of the HN-stratum of $\cG^\prime_k$ intersects the HN-stratum of $\cG_k$. 
\end{kor}
\begin{proof}
By the last remark, it is sufficient to show that $\cG_L=\cQ_k/U\times^L G$ is contained in the closure of $\cG_k^\prime$ in $\Bun_G$.

But by the definition of our cocycle, $\cG^\prime_k$ has a reduction $\cQ^-_k$ to $Q^-$ such that $\cQ^-_k/U^-\cong \cL$.
Applying the previous lemma to the $Q^-$-bundle $\cQ^-_{\lambda}$ defined above, we get a family of $G$-bundles on $\bA^1$, such that the restriction to $\bG_m$  is isomorphic to $\cG^\prime_k \times \bG_m$ and the fiber over $0$ is isomorphic to $\cQ^\prime_k/{U^-}\times^L G = \cG_L$.
\end{proof}

Now by Lemma \ref{cocycle} (4) we know that $\cQ^-_k$ has a canonical reduction to a $LU_i^-$-bundle $\cQ^-_{i,k}$. 

Let $\mu_j:=\ker(Z(L)^\circ \to \Aut(LU_i^-))$ be the kernel 
of the conjugation action and let $\overline{\lambda}:\bG_m\cong Z(L)^\circ/\mu_j\to \Aut(LU_i^-)$ be the induced map, which again extends to a map $\widetilde{\overline{\lambda}}:\bA^1\to \End(LU_i^-)$.

As before we define $\cQ^-_{i,\overline{\lambda}}:=(\cQ_{i,k}\times \bA^1)\times^Q Q$ the deformation of $\cQ_{i,k}$ into the bundle induced from its Levi quotient. 

Let $\mu_l$ be the kernel of the action on $U_i^-/U_{i+1}^-$ induced by $\overline{\lambda}$. Note that $l=1$, unless $i\geq 2$ and $U_3^-\neq \{1\}$. As remarked before, the latter can only occur if $G$ has a factor which is an exceptional group.

\begin{lemma}\label{deformationclass}
The family $\cQ^-_{i,\overline{\lambda}}|_{\Spec(k[t]/t^{l})}$ is constant. Thus  $\cQ^-_{i,\overline{\lambda}}|_{\Spec(k[t]/t^{l+1})}$ defines an element in $H^1(C,^\cL(\cu_i^-/\cu_{i+1}^-))$.
This element is the same as the one given by $\cQ^-_i/U_{i+1}^-$ under the canonical isomorphism $\cu_i^-/\cu_{i+1}^- \to U_i^-/U_{i+1}^-$ defined by the \'epinglage of $G$.
\end{lemma}
\begin{proof}
This follows immediately from the description of $\cQ^-_{i,\overline{\lambda}}$ by cocycles.
\end{proof}

{\em Now we apply Lemma \ref{adjoint} and assume that $G$ is adjoint and either contains no factor that is an exceptional group, or that $G$ is an exceptional simple group.}

\begin{kor}
$\cG^\prime_k$ is not contained in the same HN-stratum as $\cG_k$.
\end{kor}
\begin{proof}
$\cG^\prime_{\overline{\lambda}}:=(\cQ_{i,\overline{\lambda}}^-\times^{LU_i^-} G)$ is a family of $G$-bundles on $C\times \bA^1$, which defines the trivial  deformation of $\cG_L$ on $\Spec k[t]/t^{l}$. Thus, the restriction of this bundle to $C\times \Spec(k[t]/t^{l+1})$ defines an element in $H^1(C,\cL \times^G \cg)$ and by the previous lemma and \ref{cocycle}(4), we know that this element is a non-zero class in $H^1(C,^\cL(\cu_i^-/\cu_{i+1}^-))$, which is a direct summand of $H^1(C,\cL\times^G \cg)$. We claim that this implies that $\cG^\prime_k$ is not contained in the same HN-stratum  as $\cG_k$:

The group $H^1(C,\cL\times^L \cq)$ describes the deformations of $\cG_L$ which lift to $\Bun_Q\to \Bun_G$. Since $H^1(C,\cL\times^L\cg)=H^1(C,\cL\times^L \cq)\oplus H^1(C,\cL\times^L \cu^-)$ we know that the above deformation of $\cG_L$ does not lift to $\Bun_Q$ in such a way that it defines the constant deformation of $\cL\times^L Q$ mod $t^l$.

If $\cG^\prime_k$ was contained in the same HN-stratum as $\cG_k$, then the family $\cG^\prime_{\overline{\lambda}}$ would be contained in this stratum.
Now, if we can show that a canonical reduction exists for this family, then the whole family can be lifted to $\Bun_Q$ and we obtain a contradiction if $l=1$.
Furthermore, $l\neq 1$ can only happen if $G$ is an exceptional group, in which case we may assume Behrend's conjecture for $G$. But then the canonical reduction has no deformations, in particular it is constant mod $t^l$, which again contradicts our assumption.

Thus we only have to construct a canonical parabolic for the family $\cG^\prime_{\overline{\lambda}}$.
Unfortunately, it is not known whether the map $\Bun_P^{ss} \to \Bun_G$ is an embedding in positive characteristic, so we do this by hand:
We know that $\cG^\prime_{\overline{\lambda}}$ is a principal $G$-bundle over $\bA^1\times C$ such that the fiber over $0$ is $\cG_L$ and all other fibers are isomorphic to $\cG^\prime_k$. Now the reductions to $P$ which have the same degree as $\cP$ are parameterized by a connected, quasi-projective scheme over $\bA^1$ which has a single point in each fiber. Since a quasi-finite scheme is a disjoint union of a finite scheme and a scheme whose components do not map surjectively to the base, we see that this scheme is finite, thus projective over $\bA^1$. Therefore any section of this scheme over $\bG_m$ extends to the whole of $\bA^1$.

We claim that under our assumptions we always have a section over $\bG_m$:

First note that the family $\cG^\prime_{\overline{\lambda}}$ is constant when we restrict it to $Z(L)^\circ \tto Z(L)^\circ/\mu_{j}\cong \bG_m$ (by Lemma \ref{deformation}). Therefore there exists a canonical reduction of the family over $Z(L)^\circ$.

Now since $G$ is adjoint, $\lambda$ is an isomorphism and therefore either $j\in\{1,2\}$ and $l=1$ or $G$ is an exceptional group. If $l=j=1$ then we have found a reduction over $\bG_m$.
If $\text{char}(k)\neq 2$ and $j\leq 2$ then the canonical reduction on $Z(L)^\circ$ descends to $\bG_m$, because the map to $\bG_m$ is \'etale.  Finally $j>2$ can only happen if $G$ is an exceptional group, for which we already proved our claim above, using Behrend's conjecture.
\end{proof}

\section{Conclusion}

The above proves our Theorem \ref{semistab}. Namely we take any extension $\cG_R$ of $\cG_K$ to $C\times R$. If this is not semistable, then we can modify $\cG_R$ as above to obtain another extension $\cG^\prime_R$ such that the special fiber $\cG^\prime_k$ lies in a Harder-Narasimhan-stratum which contains $\cG_k$ in its closure, therefore it is contained in a stratum of higher dimension. Thus after finitely many modifications we obtain a semistable bundle.

\end{document}